\documentclass[12pt]{article}
\usepackage{amsmath}
\usepackage{amssymb}
\usepackage[dvips]{graphicx}
\font\tenmsb=msbm10

\def\Bbb#1{\hbox{\tenmsb#1}}

\font\tencmmib=cmmib10 \skewchar\tencmmib '60
\newfam\cmmibfam
\textfont\cmmibfam=\tencmmib

\font\tenmsb=msbm10

\def\Bbb#1{\hbox{\tenmsb#1}}

\def\bbox{\quad\hbox{\vrule \vbox{\hrule \vskip2pt \hbox{\hskip2pt
\vbox{\hsize=1pt}\hskip2pt} \vskip2pt\hrule}\vrule}}
\def\lessim{\ \lower4pt\hbox{$
\buildrel{\displaystyle <}\over\sim$}\ }
\def\gessim{\ \lower4pt\hbox{$\buildrel{\displaystyle >}
\over\sim$}\ }

\def\qed{\hfill\break\rightline{$\bbox$}}
\newtheorem{theorem}{Theorem}
\newtheorem{corollary}{Corollary}
\makeatletter
\@addtoreset{equation}{section}

\makeatother

\begin{document}

\title{
Some extensions of an inequality of Vapnik
and Chervonenkis.}

\author{
Dmitriy Panchenko\thanks{
The work was done during summer internship at AT\&T Research Labs at the Laboratory of Speech and Image Processing.}
\\
Department of Mathematics and Statistics
\\
The University of New Mexico\\
panchenk@math.unm.edu\\
http://www.math.unm.edu/\~{}panchenk/
}
\date{\bf August 2001}

\maketitle
\begin{abstract}
The inequality of Vapnik and Chervonenkis
controls the expectation of the function 
by its sample average uniformly over a VC-major class of functions
taking into account the size of the expectation.
Using Talagrand's kernel method we prove a similar result for the classes of
functions for which Dudley's uniform entropy integral or bracketing
entropy integral is finite. 

\end{abstract}

\vskip 5mm


\section{Introduction and main results.}
Let $\Omega$ be a measurable space with a probability
measure $P$ and $\Omega^n$ be a product space with a product
measure $P^n .$
Consider a family of measurable functions 
${\cal F}=\{f:\Omega\to [0,1]\}.$
Denote
$$
P f = \int f dP,\,\,\,\,\,\, 
\bar{f}=n^{-1}\sum_{i=1}^{n} f(x_i),\,\,\,\,\,\,
x=(x_1,\ldots,x_n)\in\Omega^n .
$$
The main purpose of this paper is to provide probabilistic bounds for
$P f$ in terms of $\bar{f}$ and the complexity assumptions on class
${\cal F}.$
We are trying to extend the following result of Vapnik
and Chervonenkis (\cite{VCer}).
Let ${\cal C}$ be a class of sets in $\Omega.$ Let
$$
S(n)=\max_{x\in \Omega^n} 
\Bigl|\Bigl\{
\{x_1,\ldots,x_n\}\cap C
: C\in{\cal C}\Bigr\}\Bigr|
$$
The VC dimension $d$ of class ${\cal C}$ is defined as
$$
d=\inf\{j\geq 1 : S(j)<2^j\}.
$$
${\cal C}$ is called VC if $d<\infty.$
The class of functions 
$\cal F$ is called VC-major if the class of sets
$$
{\cal C}=\Bigl\{\{x\in\Omega : f(x)\leq t\} : 
f\in{\cal F}, t\in R\Bigr\}
$$
is a VC class of sets in $\Omega,$ and the VC dimension of
$\cal F$ is defined as the VC dimension of ${\cal C}.$
The inequality of Vapnik and Chervonenkis states that 
(see Theorem 5.3 in \cite{Vap})
if ${\cal F}$ is
a VC-major class of $[0,1]$ valued functions
with dimension $d$ then for all $\delta>0$
with probability at least $1-\delta$ for all $f\in {\cal F}$
\begin{equation}
\frac{1}{n (Pf)^{1/2}}\sum_{i=1}^{n}(Pf - f(x_i))\leq
2\Bigl(
\frac{1}{n}\log S(2n) + 
\frac{1}{n}\log\frac{4}{\delta}
\Bigr)^{1/2},
\label{Vapnik}
\end{equation}
where for $n\geq d,$ $S(n)$ can be bounded by
$$
S(n)\leq \Bigl(\frac{en}{d}\Bigr)^{d}
$$
(see \cite{VaCher}) to give
\begin{equation}
\frac{1}{n(Pf)^{1/2}}\sum_{i=1}^{n}(Pf - f(x_i))\leq
2\Bigl(
\frac{d }{n}\log \frac{2en}{d} + 
\frac{1}{n}\log\frac{4}{\delta}
\Bigr)^{1/2}.
\label{Vapnik2}
\end{equation}
The factor $(Pf)^{-1/2}$ allows interpolation between 
the $n^{-1}$ rate for $Pf$ in the optimistic 
zero error case $\bar{f}=0$
and the $n^{-1/2}$ rate in the pessimistic case when $\bar{f}$ is ``large''.
In this paper we will prove a bound of a similar nature
under different assumptions on the complexity of the class ${\cal F}.$
Using Talagrand's abstract concentration inequality in product spaces
and the related kernel method for empirical processes \cite{Ta}
we will first prove a general result that interpolates between
optimistic and pessimistic cases.
Then we will give examples of application of this general result
in two situations when it is assumed that either Dudley's uniform
entropy integral is finite or the bracketing entropy integral is finite.

Let us formulate Talagrand's concentration inequality
that is used in the proof of our main Theorem 2 below. 
Consider a probability measure
$\nu$ on $\Omega^n$ and $x\in \Omega^n.$ 
We will denote by $x_i$ the $i^{\mbox{\scriptsize th}}$
coordinate of $x.$
If ${\cal C}_i=\{y\in\Omega^n : y_i\not= x_i\},$ 
we consider the image of the restriction of $\nu$
to ${\cal C}_i$ by the map $y\to y_i,$ and its
Radon-Nikodym derivative $d_i$ with respect to $P.$
As in \cite{Ta} we  assume that $\Omega$ is finite and each
point is measurable with a positive measure.
Let $m$ be a number of atoms in $\Omega$ and
$p_1,\ldots,p_m$ be their probabilities.
By the definition of $d_i$ we have
$$
\int_{{\cal C}_i}g(y_i)d\nu(y)=
\int_{\Omega}g(y_i)d_i(y_i)dP(y_i).
$$
For $\alpha>0$ we define a function 
$\psi_{\alpha}(x)$ by
\[
\psi_{\alpha}(x)=\left\{
\begin{array}{cl}
x^2/(4\alpha), & \mbox{when $x\leq 2\alpha,$} \\
x-\alpha, & \mbox{when $x\geq 2\alpha.$}
\end{array}
\right.
\]
We set
$$
m_{\alpha}(\nu,x)=\sum_{i\leq n}\int\psi_{\alpha}(d_i)dP\,\,\,
\mbox{ and }\,\,\,
m_{\alpha}(A,x)=\inf\{m_{\alpha}(\nu,x) : \nu(A)=1\}.
$$
For each $\alpha>0$ let $L_{\alpha}$ be any positive number
satisfying the following inequality:
\begin{equation}
\frac{2L_{\alpha}(e^{1/L_{\alpha} }-1)}{1+2L_{\alpha}}\leq\alpha.
\label{Lalpha}
\end{equation}
The following theorem holds (see \cite{Pa}).
\begin{theorem}
Let $\alpha>0$ and $L_{\alpha}$ satisfy (\ref{Lalpha}). 
Then for any $n$ and $A\subseteq \Omega^n$
we have
\begin{equation}
\int\exp
\frac{1}{L_{\alpha}}m_{\alpha}(A,x)
dP^n (x)\leq \frac{1}{P^n(A)}.
\end{equation}
\end{theorem}
Below we will only use this theorem for $\alpha=1$ and $L_1\approx 1.12.$
Let us introduce the normalized empirical process as
$$
Z(x)=
\sup_{\cal F}
\frac{1}{\varphi(f)}
\sum_{i=1}^{n} (P f -f(x_i))
,\,\,\,\,\,
x\in\Omega^n,
$$
where $\varphi : {\cal F}\to (0,\infty)$ is a function such that
$Z$ has a finite median 
$M=M(Z)<\infty,$ i.e.
\begin{equation}
P(Z\geq M)\leq \frac{1}{2}\,\,\,
\mbox{ and }\,\,\,
\forall \varepsilon>0 \,\,\,\,\,\,
P(Z\geq M+\varepsilon)<\frac{1}{2}. 
\label{median}
\end{equation}
The factor $\varphi(f)$ will play the same role as
$(n Pf)^{1/2}$ plays in (\ref{Vapnik})
The following theorem holds.

\begin{theorem}
Let $L\approx 1.12.$
If (\ref{median}) holds 
then for any $u>0,$
\begin{equation}
{\Bbb P}\Bigl(
\exists f\in{\cal F}\,\,
\sum_{i\leq n} (P f -f(x_i)) \geq M\varphi(f) + 
2\sqrt{L n u P f }
\Bigr)\leq
2e^{-u}
\label{twoone}
\end{equation}

\end{theorem}
{\bf Proof.}
The proof of the theorem repeats the proof of Theorem 2
in \cite{Pa} with some minor modifications, 
but we will give it here for completeness.
Let us consider the set 
$A=\{Z(x)\leq M\}.$
Clearly, $P^n (A)\geq 1/2.$ 
Let us fix a point $x\in \Omega^n$ and 
then choose $f\in{\cal F}.$
For any point $y\in A$ we have
$$
\frac{1}{\varphi(f)}\sum_{i=1}^{n}(P f - f(y_i))\leq M.
$$
Therefore, for any probability measure $\nu$ such that
$\nu(A)=1$ we will have
\begin{eqnarray*}
&&
\frac{1}{\varphi(f)}\sum_{i\leq n}(P f - f(x_i))-
M \leq
\frac{1}{\varphi(f)}\int \Bigl(
\sum_{i\leq n}(P f - f(x_i))-
\sum_{i\leq n}(P f - f(y_i))
\Bigr)d\nu(y)
\\
&&
=
\frac{1}{\varphi(f)}\sum_{i\leq n}\int 
(f(y_i)-f(x_i))d_i(y_i)dP(y_i).
\end{eqnarray*}
It is easy to observe that for $v\geq 0,$ and
$-1\leq u\leq 1,$
\begin{equation}
uv\leq u^2I(u> 0)+\psi_1(v).
\label{convconj}
\end{equation}
Therefore, for any $\delta> 1$
\begin{eqnarray*}
&&
\sum_{i\leq n}(P f - f(x_i))-
M\varphi(f) \leq
\delta \sum_{i\leq n}\int 
\frac{f(y_i)-f(x_i)}{\delta}d_i(y_i)dP(y_i)
\\
&&
\leq
\frac{1}{\delta}
\sum_{i\leq n}\int (f(y_i)-f(x_i))^{2}I(f(y_i)> f(x_i))dP(y_i) + 
\delta \sum_{i\leq n}\int \psi_1(d_i)dP
\end{eqnarray*}
Taking the infimum over $\nu$ we obtain that for any 
$\delta>1$
$$
\sum_{i\leq n}(P f - f (x_i)) \leq 
M\varphi(f) + 
\frac{1}{\delta} 
\sum_{i\leq n}\int (f(y_i)-f(x_i))^{2}I(f(y_i)> f(x_i))dP(y_i)
+\delta m_1(A,x). 
$$
Let us denote the random variable $\xi=f(y_1),$ 
$F_{\xi}(t)$ - the distribution function of $\xi,$ 
and $c_i=f(x_i).$ For $c\in[0,1]$ define the function $h(c)$ as
$$
h(c)=\int (f(y_1)-c)^{2}I(f(y_1)> c)dP(y_1)=
\int_{c}^{1}(t-c)^2 dF_{\xi}(t).
$$
One can check that $h(c)$ is decreasing, convex, 
$h(0)=P f^2$ and $h(1)=0.$ Therefore,
$$
\frac{1}{n}\sum_{i\leq n}h(c_i)\leq 
\Bigl(\frac{1}{n}\sum_{i\leq n}c_i\Bigr) h(1)+
\Bigl(1- \frac{1}{n}\sum_{i\leq n}c_i\Bigr) h(0)=
(1- \bar{f}) P f^2.
$$
Hence, we showed that
$$
\sum_{i\leq n}(P f - f(x_i)
\leq M\varphi(f) 
+\frac{1}{\delta}n P f + \delta m_1(A,x).
$$
Theorem 1 then implies via the application
of Chebyshev's inequality that with probability at least
$1-2 e^{-u},$
$m_1(A,x)\leq Lu$ and, hence 
$$
\sum_{i\leq n}(P f - f(x_i)
\leq M\varphi(f) + \inf_{\delta>1}
\Bigl(\frac{1}{\delta}n P f +\delta Lu\Bigr).
$$
For $u\leq nPf /L$ the infimum over $\delta>1$ equals
$2\sqrt{L n u Pf}.$ On the other hand, for $ u\geq nPf /L$
this infimum is greater than $2n Pf$ whereas the left-hand 
side is always less than $n Pf.$ 

\qed

We will now give two examples of normalization $\varphi(f)$
where we can prove that (\ref{median}) holds.

\subsection{ Uniform entropy conditions.}
Given a probability distribution $Q$ on $\Omega$
we denote
$$
d_{Q,2}(f,g)=(Q(f-g)^2)^{1/2}
$$
an $L_2-$distance on ${\cal F}$ with respect to $Q.$
Given $u>0$ we say that a subset ${\cal F}'\subset {\cal F}$ 
is $u-$separated if for any $f\not =g\in {\cal F}'$ we have
$ d_{Q,2}(f,g)> u.$ Let the {\it packing number} 
$D({\cal F},u,L_2(Q))$ be the maximal cardinality of 
any $u-$separated set. 
We will say that ${\cal F}$ satisfies the uniform entropy condition if 
\begin{equation}
\int_{0}^{\infty}\sqrt{\log D({\cal F},u)} du <\infty 
\label{Dudley},
\end{equation}
where
$$
\sup_{Q}D({\cal F},u,L_2(Q))\leq D({\cal F},u)
$$
and the supremum is taken over all discrete probability measures.
It is well known (see, for example, \cite{Du})
that if one considers the subset 
${\cal F}_p=\{f\in{\cal F} : P f\leq p\},$ then the
expectation of $\sup_{{\cal F}_p}\sum(P f - f(x_i))$
can be estimated (in some sense, since the symmetrization 
argument is required) by
\begin{equation}
\varphi(p)=\sqrt{n} \int\limits_{0}^{\sqrt{p}}\sqrt{\log D({\cal F},u)} du.
\label{phi}
\end{equation}
We will prove that it holds for all $p>0$ simultaneously.

\begin{theorem}
Assume that $D({\cal F},1)\geq 2$
and (\ref{Dudley}) holds. 
If $\varphi$ is defined by (\ref{phi})
then the median 
$$
M=M\Bigl(
\sup_{\cal F}
\frac{1}{\varphi(P f)}
\sum_{i=1}^{n} (P f -f(x_i))
\Bigr)
\leq K<\infty
,
$$
is finite, where $K$ is an absolute constant.

\end{theorem}
{\bf Proof.} The proof is based on standard symmetrization
and chaining techniques. 
We will first prove that
\begin{equation}
{\Bbb P}\Bigl(\sup_{\cal F}\frac{\sum (P f -f(x_i))}
{\varphi(P f)}\geq u\Bigr)\leq
2 {\Bbb P}\Bigl(\sup_{\cal F}\frac{\sum (f(y_i) -f(x_i))}
{\varphi(\bar{f}(x,y))}
\geq u-\Bigl(\frac{2}{\log 2}\Bigr)^{1/2}\Bigr).
\label{symm}
\end{equation}
where
$$
\bar{f}(x,y)=\frac{1}{2n}\sum(f(y_i)+f(x_i)).
$$
Let 
$$
A=\Bigl\{x :\sup_{\cal F}\frac{\sum (P f -f(x_i))}
{\varphi(P f)}\geq u \Bigr\}.
$$
Let $x\in A$ and $f\in{\cal F}$ be such that
$\sum (P f -f(x_i))/\varphi(P f)\geq u.$
Chebyshev's inequality implies 
$$
{\Bbb P}\Bigl(
|\sum_{i=1}^{n}(Pf -f(y_i))|\geq \sqrt{2n Pf }
\Bigr)\leq
\frac{n\mbox{Var}f}{2 n Pf}\leq \frac{1}{2},
$$
where $y=(y_1,\ldots,y_n)$ lives on an independent 
copy of $(\Omega^n,P^n).$ We will show that the inequalities
$$
n Pf\leq \sum f(y_i) +\sqrt{2n Pf },\,\,\,\,\,\,
u\leq \frac{\sum (P f -f(x_i))}{\varphi(P f)}
$$
imply that
$$
\frac{\sum (f(y_i) -f(x_i))}
{\varphi(\bar{f}(x,y))}
\geq u-\Bigl(\frac{2}{\log 2}\Bigr)^{1/2} .
$$
If we define by ${\Bbb P}_y$ the probability measure
on the space of $y,$ it would mean that 
\begin{eqnarray*}
&&
\frac{1}{2} I(x\in A)\leq
{\Bbb P}_{y}\Bigl(
|\sum_{i=1}^{n}(Pf -f(y_i))|\geq \sqrt{2n Pf}
\Bigr)
\leq
{\Bbb P}_{y}\Bigl( 
\frac{\sum (f(y_i) -f(x_i))}
{\varphi(\bar{f}(x,y))}
\geq u-\Bigl(\frac{2}{\log 2}\Bigr)^{1/2}
\Bigr)
\\
&&
\leq
{\Bbb P}_{y}\Bigl( 
\sup_{\cal F}\frac{\sum (f(y_i) -f(x_i))}
{\varphi(\bar{f}(x,y))}
\geq u-\Bigl(\frac{2}{\log 2}\Bigr)^{1/2}
\Bigr)
\end{eqnarray*}
and taking expectation of both sides 
with respect to $x$
would prove (\ref{symm}).
To show the remaining implication we consider two cases when
$nPf\leq\sum f(y_i)$ and $nPf\geq\sum f(y_i).$
First assume that $nPf\leq\sum f(y_i).$
Since, as easily checked, both $\varphi(p)$ and
$p/\varphi(p)$ are increasing we get
$$
\frac{\sum (Pf-f(x_i))}{\varphi(Pf)}\leq
\frac{\sum (f(y_i)-f(x_i))}{\varphi(n^{-1}\sum f(y_i))}\leq
\frac{\sum (f(y_i)-f(x_i))}{\varphi(\bar{f}(x,y))}.
$$
In the case $nPf\geq\sum f(y_i)$ we have
$$
\frac{\sum (Pf-f(x_i))}{\varphi(Pf)}\leq
\frac{\sum (f(y_i)-f(x_i))}{\varphi(Pf)}+
\frac{\sqrt{2n Pf}}{\varphi(Pf)}\leq
\frac{\sum (f(y_i)-f(x_i))}{\varphi(\bar{f}(x,y))}
+\frac{\sqrt{2n Pf}}{\varphi(Pf)}.
$$
The assumption
$D({\cal F},\sqrt{Pf})\geq D({\cal F},1)\geq 2$
garantees that
$\varphi(Pf)\geq \sqrt{n Pf \log 2}$
and, finally,
$$
u\leq
\frac{\sum (f(y_i)-f(x_i))}{\varphi(\bar{f}(x,y))}
+\Bigl(\frac{2}{\log 2}\Bigr)^{1/2}
$$
which completes the proof of (\ref{symm}).
We have
$$
{\Bbb P}\Bigl(\sup_{\cal F}\frac{\sum ( f(y_i) -f(x_i))}
{\varphi(\bar{f}(x,y))}\geq u\Bigr)=
{\Bbb E}
{\Bbb P}_{\varepsilon}
\Bigl(\sup_{\cal F}\frac{\sum \varepsilon_i( f(y_i) -f(x_i))}
{\varphi(\bar{f}(x,y))}\geq u\Bigr),
$$
where $(\varepsilon_i)$ is a sequence of Rademacher random variables.
We will show that there exists $u$ independent of $n$ such that for any 
$x,y\in\Omega^n$
$$
{\Bbb P}_{\varepsilon}
\Bigl(\sup_{\cal F}\frac{\sum \varepsilon_i( f(y_i) -f(x_i))}
{\varphi(\bar{f}(x,y))}\geq u\Bigr)< \frac{1}{2}.
$$
Clearly, this will prove the statement of the theorem.
For a fixed $x,y\in\Omega^n$ let
$$
F=\{(f(x_1),\ldots,f(x_n),f(y_1),\ldots,f(y_n)) : f\in{\cal F}\}\subset R^{2n}
$$
and
$$
d(f,g)=\Bigl(\frac{1}{2n}\sum_{i=1}^{2n}(f_i-g_i)^2\Bigr)^{1/2},
\,\,\,\,\,\,\,\,\,
f,g\in F.
$$
The packing number of $F$ with respect to $d$ can be bounded by
$D(F,u,d)\leq D({\cal F},u).$
Consider an increasing sequence of sets
$$
\{0\}=F_0\subseteq F_1\subseteq F_2\subseteq \ldots
$$
such that for any $g\not =h\in F_j,$ $d(g,h)>2^{-j}$
and for all $f\in F$ there exists $g\in F_j$ such that
$d(f,g)\leq 2^{-j}.$ The cardinality of $F_j$ can be bounded by
$$
|F_j|\leq D(F,2^{-j},d)\leq D({\cal F},2^{-j}).
$$
For simplicity of notations we will write $D(u):= D({\cal F},u).$
If $D(2^{-j})=D(2^{-j-1})$ then in the construction of the sequence
$(F_j)$ we will set $F_j$ equal to $F_{j+1}.$
We will now define the sequence of projections 
$\pi_j:F\to F_j,\,j\geq 0$ in the following way.
If $f\in F$ is such that $d(f,0)\in(2^{-j-1},2^{-j}]$
then set $\pi_0(f)=\ldots=\pi_j(f)=0$ and for $k\geq j+1$ 
choose $\pi_k(f)\in F_k$ such that $d(f,\pi_k(f))\leq 2^{-k}.$
In the case when $F_k=F_{k+1}$ we will choose $\pi_{k}(f)=\pi_{k+1}(f).$
This construction implies that $d(\pi_{k-1}(f),\pi_{k}(f))\leq 2^{-k+2}.$
Let us introduce a sequence of sets
$$
\Delta_j=\{g-h : g\in F_j, h\in F_{j-1}, d(g,h)\leq 2^{-j+2}\},\,\,\,\,
j\geq 1,
$$
and let $\Delta_j = \{0\}$ if $D(2^{-j})=D(2^{-j+1}).$
The cardinality of $\Delta_j$ does not exceed 
$$
|\Delta_j|\leq |F_j|^2\leq D(2^{-j})^2.
$$
By construction any $f\in F$ can be represented 
as a sum of elements from $\Delta_j$
$$
f=\sum_{j\geq 1}(\pi_j(f)-\pi_{j-1}(f)),\,\,\,\,\,\,
\pi_{j}(f)-\pi_{j-1}(f) \in \Delta_j.
$$
Let
$$
I_j=\sqrt{n}\int\limits_{2^{-j-1}}^{2^{-j}}\sqrt{\log D(u)} du
$$
and define the event
$$
A=
\bigcup_{j=1}^{\infty}
\{
\sup_{f\in\Delta_j} \sum_{i=1}^{n}\varepsilon_i (f_{i+n}-f_{i})
\geq
u I_j
\}.
$$
On the complement $A^c$ of the event $A$ we have for any 
$f\in F$ such that $d(f,0)\in(2^{-j-1},2^{-j}]$
\begin{eqnarray*}
&&
\sum_{i=1}^{n}\varepsilon_i (f_{i+n}-f_{i})=
\sum_{k\geq j+1}\sum_{i=1}^{n}\varepsilon_i 
\bigl(
(\pi_k(f)-\pi_{k-1}(f))_{i+n}-(\pi_k(f)-\pi_{k-1}(f))_{i}
\bigr)
\\
&&
\leq
\sum_{k\geq j+1} u I_k \leq
u\sqrt{n}\int\limits_{0}^{2^{-j-1}}\sqrt{\log D(u)} du\leq
u\sqrt{n}\int\limits_{0}^{(\bar{f})^{1/2}}\sqrt{\log D(u)} du,
\end{eqnarray*}
where $\bar{f}=(2n)^{-1}\sum_{i\leq 2n} f_i,$ since 
$2^{-j-1}<d(f,0)\leq (\bar{f})^{1/2}.$
It remains to prove that for some absolute constant $u,$
$P(A)<1/2.$ Indeed,
\begin{eqnarray*}
&&
P(A)\leq \sum_{j=1}^{\infty}
P\bigl(
\sup_{f\in\Delta_j} \sum_{i=1}^{n}\varepsilon_i (f_{i+n}-f_{i})
\geq u I_j
\bigr) 
\\
&&
\leq
\sum_{j=1}^{\infty}
|\Delta_j|\exp \Bigl\{
-\frac{u^2 I_j^2}{n 2^{-2j+6}}
\Bigr\}
I\bigl(D(2^{-j})>D(2^{-j+1})\bigr)
\\
&&
\leq
\sum_{j=1}^{\infty}
\exp \Bigl\{
2\log D(2^{-j})
-\frac{u^2 I_j^2}{n 2^{-2j+6}}
\Bigr\}
I\bigl(D(2^{-j})>D(2^{-j+1})\bigr),
\end{eqnarray*}
since for $f\in\Delta_j$
$$
\sum_{i=1}^{n}(f_{i+n}-f_i)^2\leq 2\sum_{i=1}^{2n} f_i^2
\leq
n 4\cdot 2^{-2j+4}.
$$
The fact that $D(u)$ is decreasing implies
$$
\frac{I_j}{\sqrt{n}\ 2^{-(j+1)}}\geq \sqrt{\log D(2^{-j})}
$$
and, therefore,
\begin{eqnarray*}
&&
P(A)\leq \sum_{j=1}^{\infty}
\exp \{
-\log D(2^{-j})(u^2 2^{-8}-2)
\}
I\bigl(D(2^{-j})>D(2^{-j+1})\bigr)
\\
&&
\leq
\sum_{j=1}^{\infty}
\frac{1}{D(2^{-j})^{\alpha}}
I\bigl(D(2^{-j})>D(2^{-j+1})\bigr)
\leq
\sum_{j=2}^{\infty}
\frac{1}{j^{\alpha}}<\frac{1}{2},
\end{eqnarray*}
for $\alpha=u^2/2^8 -2$ big enough.

\qed

Combining Theorem 2 and Theorem 3 we get

\begin{corollary}
If (\ref{Dudley}) holds then there exists
an absolute constant $K>0$ such that for any $u>0$ 
with probability at least $1-2 e^{-u}$ for all
$f\in {\cal F}$
$$
\sum_{i=1}^{n}(Pf -f(x_i))\leq 
K\Bigl(
\sqrt{n}\int\limits_{0}^{(Pf)^{1/2}}\sqrt{\log D({\cal F},u)}du +
\sqrt{n u P f}
\Bigr). 
$$

\end{corollary}

\subsection{Bracketing entropy conditions.}
Given two functions $g,h:\Omega\to [0,1]$ 
such that $g\leq h$ and $(P(h-g)^2)^{1/2}\leq u$
we will  call a set of all functions $f$ such that
$g\leq f\leq h$ a $u-$bracket with respect to $L_2(P).$
The $u-$bracketing number $N_{[]}({\cal F},u, L_2(P))$ is the
minimum number of $u-$brackets needed to cover ${\cal F}.$ 
Assume that
\begin{equation}
\int\limits_{0}^{\infty} \sqrt{\log N_{[]}({\cal F},u,L_2(P))} du<\infty
\label{bracket}
\end{equation}
and denote
$$
\varphi(p)=
\sqrt{n}
\int\limits_{0}^{\sqrt{p}} \sqrt{\log N_{[]}({\cal F},u,L_2(P))} du.
$$
Then the following theorem holds.
\begin{theorem}
Assume that $N_{[]}({\cal F},1,L_2(P))\geq 2$
and (\ref{bracket}) holds. 
If $\varphi$ is defined by (\ref{phi})
then the median 
$$
M=M\Bigl(
\sup_{\cal F}
\frac{1}{\varphi(P f)}\sum_{i=1}^{n} (P f -f(x_i))
\Bigr)
\leq K({\cal F})<\infty
,
$$
where $K({\cal F})$ does not depend on $n.$
\end{theorem}

We omit the proof of this theorem since it is a modification
of a standard bracketing entropy bound (see Theorem 2.5.6 
and 2.14.2 in \cite{Well})
similar to what Theorem 3
is to the standard uniform entropy bound. The argument is more subtle
as it involves a truncation argument required by the application of 
Bernstein's inequality but otherwise it repeats Theorem 3.
Combining Theorem 2 and Theorem 4 we get 
\begin{corollary}
If (\ref{bracket}) holds then there exists
an absolute constant $K>0$ such that for any $u>0$ 
with probability at least $1-2 e^{-u}$ for all
$f\in {\cal F}$
$$
\sum_{i=1}^{n}(Pf -f(x_i))\leq 
K \Bigl(
\sqrt{n}\int\limits_{0}^{(Pf)^{1/2}}
\sqrt{\log N_{[]}({\cal F},u,L_2(P))}du +
\sqrt{ n u P f}
\Bigr). 
$$

\end{corollary}

\section{Examples of application.}

{\bf Example 1} (VC-subgraph classes of functions).
A class of functions ${\cal F}$ is called VC-subgraph
if the class of sets 
$$
{\cal C}=
\Bigl\{
\{(\omega,t): 
\omega\in\Omega, t\in R, t\leq f(\omega)\} : f\in {\cal F}
\Bigr\}
$$
is a VC-class of sets in $\Omega\times R .$ 
The VC dimension of $\cal F$ is equal to the VC dimension $d$ of $\cal C.$
On can use Corollary 3 in \cite{Haussler} to show that
$$
D({\cal F},u) \leq e(d+1)\Bigl(\frac{2e}{u^2}\Bigr)^d.
$$
Corollary 1 implies in this case that for any $\delta>0$
with probability at least $1-\delta$ for all $f\in{\cal F}$
\begin{equation}
\frac{1}{n(Pf)^{1/2}}\sum_{i=1}^{n}(Pf-f(x_i))\leq
K\Bigl(
\Bigl(
\frac{d }{n}\log n
\Bigr)^{1/2}
+
\Bigl(\frac{1}{n}\log\frac{1}{\delta}\Bigr)^{1/2}
\Bigr),
\label{subgraph}
\end{equation}
where $K>0$ is an absolute constant. 
Instead of the $\log n$ on the right-hand side of (\ref{subgraph})
one could also write $\log(1/Pf),$ but we simplify the bound
to eliminate this dependence on $Pf.$ 
Note that the bound is similar to the bound (\ref{Vapnik2})
for VC classes of set and VC-major classes.
Unfortunately, our proof does not allow us to recover 
the same small value of $K=2$ as for VC classes of sets.

(\ref{subgraph}) 
improves the main result in \cite{Li}, where it was shown that
for any fixed $\nu>0$ for any $\delta>0$ with probability at least
$1-\delta$ for all $f\in {\cal F}$
\begin{equation}
\frac{\sum(Pf - f(x_i))}{\sum(Pf + f(x_i)) +n\nu}\leq
K\Bigl(
\frac{1}{n\nu}\Bigl(
d\log\frac{1}{\nu}+\log\frac{1}{\delta}
\Bigr)
\Bigr)^{1/2}.
\label{li}
\end{equation}
It is easy to see that, in a sense, one would get (\ref{subgraph})
from (\ref{li}) only after optimizing over $\nu.$
Indeed, for $Pf\lessim \nu,$ (\ref{li}) gives
$$
\frac{1}{n}\sum(Pf - f(x_i))\lessim 
\Bigl(
\frac{\nu}{n}\Bigl(
d\log\frac{1}{\nu}+\log\frac{1}{\delta}
\Bigr)
\Bigr)^{1/2},
$$
which is implied by (\ref{subgraph}) as well. 
For $Pf \gtrsim \nu,$
(\ref{li}) gives 
$$
\frac{1}{n}\sum(Pf - f(x_i))\lessim 
\Bigl(\frac{Pf}{\nu}\Bigr)^{1/2}
\Bigl(
\frac{Pf}{n}\Bigl(
d\log\frac{1}{\nu}+\log\frac{1}{\delta}
\Bigr)
\Bigr)^{1/2},
$$
which compared to (\ref{subgraph}) contains an additional factor
of $(Pf/\nu)^{1/2}.$ In the situation when $\nu$ is small 
(this is the only interesting case) this factor introduces
an unnecessary penalty for any function $f$ such that
$Pf\gg \nu.$ Hence, for a fixed $\nu$ (\ref{li}) improves
the bound for $Pf\leq\nu$ at cost of $f$ with $Pf\geq\nu.$ 

One can find alternative extensions of
(\ref{li}) in \cite{Kohler}. 
For some other applications of Corollary 1 see \cite{Pa2}.

{\bf Example 2} (Bracketing entropy). Assume that either
$$
D({\cal F},u)\leq cu^{-\gamma} 
\mbox{ or } 
N_{[]}({\cal F},u,L_2(P))\leq cu^{-\gamma},\,\,\,
\gamma\in(0,2). 
$$
Then Corollary 1 or Corollary 2 imply that for $u>0$
with probability at least $1-2e^{-u}$ for all $f\in{\cal F}$
$$
Pf-\bar{f}\leq\frac{c_{\gamma}}{\sqrt{n}}
\bigl(
(Pf)^{\frac{1}{2}-\frac{\gamma}{4}}
+
(uPf)^{\frac{1}{2}}
\bigr).
$$
If $\bar{f}=0$ then it is easy to see that for
$u\leq n^{\frac{\gamma}{2+\gamma}}$ we have 
$$
Pf\leq K_{\gamma} n^{-\frac{2}{2+\gamma}}.
$$
As an example, if ${\cal F}$ is a class of indicator functions
for sets with $\alpha-$smooth boundary in $[0,1]^l$ 
and $P$ is Lebesgue absolutely continuous with bounded density
then well known bounds on the bracketing entropy 
due to Dudley (see \cite{Du}) imply that $\gamma=2(l-1)/\alpha$
and $Pf\leq K_{\alpha}n^{-\frac{\alpha}{l-1+\alpha}}.$
Even though $\gamma=2(l-1)/\alpha$ may be greater than $2$
and Corollary 2 is not immediately applicable, one can 
generalize Theorem 4 to different choices of $\varphi(x),$
using the standard truncation in the chaining argument, 
to obtain the above rates even for $\gamma\geq 2.$

{\bf Acknowledgments.} We would like to thank the referee for 
several very helpful comments and suggestions.

\end{document}